\documentclass[a4paper, 12pt]{article}

\usepackage[koi8-r,cp1251]{inputenc}
\usepackage{amsfonts,amssymb,amsmath, hyperref}
\usepackage[final]{epsfig}
\usepackage{graphicx}

\begin{document}

\begin{center}
\textbf{ HAUSDORFF OPERATORS ON SOME SPACES OF HOLOMORPHIC FUNCTIONS ON THE UNIT DISC}
\end{center}

\

\begin{center}
\textbf{A. R. Mirotin}
\end{center}

\begin{center}
 amirotin@yandex.ru
\end{center}

\

\textbf{Abstract.} {\small We introduce  Hausdorff operators over the unit disc and give conditions for boundedness of such operator in Bloch, Bergman,  and Hardy spaces on the disc. Identity approximation by Hausdorff operators is also considered.}

\


Mathematics Subject Classification 47G10, 47B38, 46E30

\

\section{Introduction and preliminaries}

In \cite{JMAA} Hausdorff operators were defined over locally compact groups $G$  via their automorphism groups $\mathrm{Aut}(G)$. But for the circle group $\mathbb{T}$ this definition  leads to almost trivial Hausdorff operators  due to the almost triviality of $\mathrm{Aut}(\mathbb{T})$.
The main idea of this note is as follows. If we want to fix the aforementioned problem we should consider  spaces of functions on the boundary $\mathbb{T}=\partial \mathbb{D}$ of the unit disc $\mathbb{D}$ with holomorphic extension into the  disc
and employ the reach  group $\mathrm{Aut}_0(\mathbb{D})$ of involutive  M\"{o}bius automorphisms  of the disc instead of $\mathrm{Aut}(\partial \mathbb{D})$.

The one-dimensional Hausdorff operators were introduced by Hardy \cite[Chapter XI]{H}. The modern theory of Hausdorff operators was inspired by \cite{LM}; see surveys \cite{Ls}, \cite{CFW}. Hardy introduced his operators as a continuous analog of Hausdorff means; similarly, to get   a mean of a holomorphic function $f$ on $\mathbb{D}$ it is natural to integrate $f$ over all natural holomorphic  mixing of the disc, i.~e., over the  elements of the group $\mathrm{Aut}(\mathbb{D})$ of M\"{o}bius transformations of $\mathbb{D}$ (see theorems 4 and 5 below). So, we arrive at the following

\textbf{Definition 1.} For a holomorphic function $f$ on $\mathbb{D}$  the value of a Hausdorff operator at $f$ is as follows:
$$
(\mathcal{H}_{K,\mu} f)(z):=\int_{\mathbb{D}}K(w)f(\varphi_w(z))d\mu(w), \ z\in \mathbb{D},
$$
where $\mu$ is some fixed positive Radon measure on $\mathbb{D}$,  $K$ is some fixed $\mu$-measurable function on $\mathbb{D}$, and $\varphi_w$ is an element of $\mathrm{Aut}(\mathbb{D})$ of the form
$$
\varphi_w(z)=\frac{w-z}{1-\overline{w}z}\quad (w\in \mathbb{D}).
$$
We shall denote the set of such elements by $\mathrm{Aut}_0(\mathbb{D})$.

\textbf{Example}. If $\mu$ is concentrated on a subsequence $\{w_n\}\subset \mathbb{D}$ we get a class of operators  of the form
$$
(\mathcal{H}_{d} f)(z):=\sum_{n=0}^\infty d_nf\left(\frac{w_n-z}{1-\overline{w_n}z}\right)
$$
(discrete Hausdorff operators over $\mathbb{D}$). If the sequence $\{w_n\}$ is finite this is a sort of so-called functional operators (see \cite{AA}).

\textbf{Remark}. In \cite{KSZhu}   operators in $L^p(dH)$ of the form
$$
\mathbb{K}f(z)=\int_{\mathbb{D}}K(w)f(\varphi_z(w))dH(w)
$$
where $dH$  is the M\"{o}bius invariant area measure  were introduced  under the name of Hausdorff-Berezin operators. This class of operators is  different from our as the previous example showes.

It should be noted that different  operators of Hausdorff type on spaces of holomorphic functions in the disc or the half-plane were considered in \cite{AL}, \cite{GP},  \cite{Stil}, \cite{ArvSis}, \cite{HKQ}, \cite{BBMM}, \cite{KS} but all of them are not special cases of the  definition 1.

\section{Boundedness of Hausdorff operators}

The following  lemma is  crucial for the  proofs of our main results.

\textbf{Lemma 1.} \cite[Lemma 2]{JMAA} \textit{Let $(X;m)$ be a measure space,
$\mathcal{F}(X)$
 some Banach space of $m$-measurable functions on $X$,  $(\Omega,\mu)$ a $\sigma$-compact quasi-metric space with  positive Radon measure $\mu$, and $F(w, z)$ a function on $\Omega\times X$. Assume that}

(a) \textit{the convergence of a sequence in norm in $\mathcal{F}(X)$ yields the convergence of some
subsequence to the same function for $m$-a.~e. $z\in X$; }

(b)  \textit{$F(w, \cdot) \in \mathcal{F}(X)$ for $\mu$-a.~e. $w\in \Omega$;}

 (c) \textit{the map  $w\mapsto F(w, \cdot):\Omega\to \mathcal{F}(X)$ is Bochner integrable with respect to } $\mu$.

  \textit{Then for $m$-a.~e. $z\in  X$ one has}
$$
\left((B)\int_\Omega F(w, \cdot)d\mu(w)\right)(z)=
\int_\Omega F(w, z)d\mu(w).
$$

\subsection{The Bloch space}

The Bloch space $\mathcal{B}$ of $\mathbb{D}$ is defined to be the space of analytic functions
$f$ on $\mathbb{D}$ such that
$$
\|f\|_{\mathcal{B}}:= \sup\{(1 - |z|^2)|f'(z)|: z \in\mathbb{ D}\} < +\infty.
$$
Then $\|\cdot\|_{\mathcal{B}}$ is a  complete semi-norm on $\mathcal{B}$, and $\mathcal{B}$ can be
made into a Banach space by introducing the norm
$$
\|f\|'_{\mathcal{B}}:= \|f\|_{\mathcal{B}}+|f(0)|.
$$
It is important that $\|\cdot\|_{\mathcal{B}}$ is M\"{o}bius invariant, that is,
$$
\|f\circ \varphi_w\|_{\mathcal{B}}= \|f\|_{\mathcal{B}}
$$
for all $\varphi_w\in \mathrm{Aut}_0(\mathbb{D})$ \cite[p. 101]{Zhu}.

\textbf{Theorem 1.} \textit{If $K(w)$ and  $K(w)\log(1-|w|)$ belong to $L^1(\mu)$, then the operator $\mathcal{H}_{K,\mu}$ is bounded on $\mathcal{B}$  and }
$$
\|\mathcal{H}_{K,\mu}\|_{\mathcal{B}\to\mathcal{ B}}\le \int_{\mathbb{D}}|K(w)|\left(2+ \frac{1}{2}\log\frac{1+|w|}{1-|w|}\right)d\mu(w).
$$

Proof. Note that the conditions of Lemma 1 hold for $X=\Omega=\mathbb{D}$, $m=A$, where $dA(z)=dxdy/\pi$ stands for the normalized area (Lebesgue) measure on $\mathbb{D}$, $\mathcal{F}(X)=\mathcal{B}$, and  $F(w,z)=K(w)f(\varphi_w(z))$. Indeed,

(a) if $f_n\in \mathcal{B}$ and $f_n\to 0$ strongly, then $f_n(0)\to 0$ and $f_n'\to 0$ uniformly on $\{z\in\mathbb{D}:|z|<r\}$ for every $r\in(0,1)$ and therefore for all $z\in \mathbb{D}$
$$
f_n(z)=\int_0^zf'_n(t)dt+f_n(0)\to 0\quad (n\to\infty);
$$

(b) since $\|f\circ \varphi_w\|_{\mathcal{B}}= \|f\|_{\mathcal{B}}$, this is obvious;

(c) this follows from the estimate
$$
\| F(w, \cdot)\|'_\mathcal{B}\le |K(w)|
\left(2+ \frac{1}{2}\log\frac{1+|w|}{1-|w|}\right)
$$
(see below).

 So, by Lemma 1 for $dA$-a.e. $z\in \mathbb{D}$
$$
\mathcal{H}_{K,\mu} f(z)=\left((B)\int_{\mathbb{D}}K(w)f\circ \varphi_wd\mu(w)\right)(z),\eqno(1)
$$
the Bochner integral for $\mathcal{B}$.  The right-hand side of this equality is continuous. To shaw that the left-hand side  of (1) is continuous, too, we choose an arbitrary point  $z_0\in \mathbb{D}$ and a compact neighborhood $U\subset \mathbb{D}$ of $z_0$. By \cite[Theorem 5.5]{Zhu} we have
$$
|f(w)|\le |f(w)-f(0)|+|f(0)|\le \|f\|_{\mathcal{B}}\beta(0,w)+|f(0)|,
$$
where $\beta(z,w)$ stands for the Bergman distance. Since $\beta$ is M\"{o}bius invariant (see, e.g.,\cite{Zhu}), it follows that  one can find such $z_1, z_2\in K$ that  for all $z\in U$
$$
|f(\varphi_w(z))|\le  \|f\|_{\mathcal{B}}\beta(0,\varphi_w(z))+|f(-z)|= \|f\|_{\mathcal{B}}\beta(\varphi_w(w),\varphi_w(z))+|f(-z)|=
$$
$$
\|f\|_{\mathcal{B}}\beta(w,z)+|f(-z)|\le \|f\|_{\mathcal{B}}(\beta(w,0)+\beta(0,z))+|f(-z)|\le
$$
$$
 \|f\|_{\mathcal{B}}\left(\frac{1}{2}\log\frac{1+|w|}{1-|w|}+\beta(0,z_1)\right)+|f(-z_2)|.
$$
Thus, if $z\in U$ we get    a $\mu$ integrable majorant for the left-hand side of (1) of the form
$$
K(w)\left(C_1\log\frac{1+|w|}{1-|w|}+C_2\right).
$$

Since both sides of the  equality (1) are continuous, it is valid for all $z\in \mathbb{D}$, i.e.,
$$
\mathcal{H}_{K,\mu} f=(B)\int_{\mathbb{D}}K(w)f\circ \varphi_wd\mu(w).
$$
Since $\|f\circ \varphi_w\|_\mathcal{B}=\|f\|_\mathcal{B}$ \cite[p. 101]{Zhu}, and $f\circ \varphi_w(0)=f(w)$, it follows that
$$
\|\mathcal{H}_{K,\mu} f\|'_{\mathcal{B}}\le\int_{\mathbb{D}}|K(w)|\|f\circ \varphi_w\|'_{\mathcal{B}}d\mu(w)=
$$
$$
\|f\|_{\mathcal{B}}\int_{\mathbb{D}}|K(w)|d\mu(w)
+\int_{\mathbb{D}}|K(w)||f(w)|d\mu(w).
$$
As was mentioned above,
$$
|f(w)|\le \|f\|_{\mathcal{B}}\beta(0,w)+|f(0)|.
$$
Therefore
$$
\|\mathcal{H}_{K,\mu} f\|'_{\mathcal{B}}\le\int_{\mathbb{D}}|K(w)|d\mu(w)\|f\|_{\mathcal{B}}+
$$
$$
\int_{\mathbb{D}}|K(w)|\beta(0,w)d\mu(w)\|f\|_{\mathcal{B}}+\int_{\mathbb{D}}|K(w)|d\mu(w)|f(0)|
\le
$$
$$
\int_{\mathbb{D}}|K(w)|(2+\beta(0,w))d\mu(w)\|f\|'_{\mathcal{B}}.
$$

This completes the proof.

\subsection{Bergman spaces}

Let $\alpha >-1$ and $dA_\alpha(z):=(\alpha+1)(1-|z|^2)^{\alpha}dA(z)$.  For $p > 0$ the  Bergman spaces with standard weights are defined by
$$
L^p_a(dA_\alpha):=H(\mathbb{D})\cap L^p(dA_\alpha)
$$
where $H(\mathbb{D})$ is the space of analytic functions in $\mathbb{D}$ (see, e.g., \cite[Chapter 4]{Zhu}).

In this subsection we shall consider the case $d\mu=dA, \|\cdot\|_{p,\alpha}$ denotes the norm in $L^p_a(dA_\alpha)$ induced from $L^p(dA_\alpha)$.

\textbf{Theorem 2.}
 \textit{Let $p\ge 1$ and the function  $K(w)/(1-|w|)^{\frac{2+\alpha}{p}}$ belongs to $L^1(dA)$.  Then the operator $\mathcal{H}_{K,A}$ is bounded on $L^p_a(dA_\alpha)$ and}
$$
\|\mathcal{H}_{K,A_\alpha}\|_{L^p_a\to L^p_a}\le \int_{\mathbb{D}}|K(w)|\left(\frac{1+|w|}{1-|w|}\right)^{\frac{2+\alpha}{p}}dA(w).
$$

Proof. According to \cite[Proposition 4.3]{Zhu},
$$
\int_{\mathbb{D}}g\circ \varphi_w(z)dA_\alpha(z)=\int_{\mathbb{D}}g(z)\frac{(1-|w|^2)^{2+\alpha}}{(1-\overline{w}z)^{2(2+\alpha)}}dA_\alpha(z)
$$
for nonnegative $g$. Putting here $g=|f|^p$ where $f\in L^p_a(dA_\alpha)$  we get
$$
\|f\circ \varphi_w\|^p_{p,\alpha}=\int_{\mathbb{D}}|f\circ \varphi_w(z)|^pdA_\alpha(z)=\int_{\mathbb{D}}|f(z)|^p\frac{(1-|w|^2)^{2+\alpha}}{|1-\overline{w}z|^{2(2+\alpha)}}dA_\alpha(z).
$$
Since $|1-\overline{w}z|\ge 1-|w|$, this implies that
$$
\|f\circ \varphi_w\|^p_{p,\alpha}\le\int_{\mathbb{D}}|f(z)|^p\frac{(1-|w|^2)^{2+\alpha}}{(1-|w|)^{2(2+\alpha)}}dA_\alpha(z)=
\left(\frac{1+|w|}{1-|w|}\right)^{2+\alpha}\|f\|^p_{p,\alpha}.\eqno(2)
$$
 Note that all the conditions of Lemma 1 are satisfied for $(X,m)=(\mathbb{D},dA_\alpha)$, $(\Omega,\mu)=(\mathbb{D},dA)$,  $\mathcal{F}(X)=L^p_a(dA_\alpha)$, and  $F(w,z)=K(w)f(\varphi_w(z))$ (indeed, (b) and (c) follow from the  estimate (2); (a) is a consequence of a well known theorem of Riesz).

  It follows in view of Lemma 1 that for $dA_\alpha$ a.e. $z\in \mathbb{D}$
 $$
\mathcal{H}_{K,\mu} f(z)=(B)\int_{\mathbb{D}}K(w)f\circ \varphi_wd\mu(w)(z),\eqno(3)
$$
the Bochner integral for $L^p_a(dA_\alpha)$. As in the proof of Theorem 1 to show that
the  equality (3) is  valid for all $z\in \mathbb{D}$ we shall prove the continuity of its left-hand side. To this end we shall prove that for every $z_0\in \mathbb{D}$ and for every compact neighborhood $U\subset \mathbb{D}$ of $z_0$ there is a $dA$ integrable majorant
for the integrand of the left-hand side of (3). First note that by \cite[Proposition 4.13]{Zhu} there is a constant $C>0$ such that for $f\in L^p_a(dA_\alpha),z\in \mathbb{D}$
$$
|f(z)|^p\le \frac{C}{(1-|z|^2)^{2+\alpha}}\|f\|^p_{p,\alpha}.
$$
It follows that
$$
|f(\varphi_w(z))|\le \frac{C^{1/p}}{(1-|\varphi_w(z)|^2)^{\frac{2+\alpha}{p}}}\|f\|_{p,\alpha}.
$$
On the other hand,
$$
1-|\varphi_w(z)|^2=\frac{(1-|w|^2)(1-|z|^2)}{|1-\overline{w}z|^2}\ge \frac{(1-|w|^2)(1-|z|^2)}{(1+|w|)^2}=
$$
$$
\frac{(1-|w|)(1-|z|^2)}{1+|w|}\ge \frac{1}{2}(1-|z|^2)(1-|w|).
$$
 Thus,
$$
|f(\varphi_w(z))|\le \frac{2^{\frac{2+\alpha}{p}}C^{1/p}}{(1-|z|^2)^{\frac{2+\alpha}{p}}(1-|w|)^{\frac{2+\alpha}{p}}}\|f\|_{p,\alpha}.
$$
This yields that for every compact neighborhood $U\subset \mathbb{D}$ of $z_0$ there is a constant $C_U$ such that
$$
|K(w)||f(\varphi_w(z))|\le \frac{C_U}{(1-|w|)^{\frac{2+\alpha}{p}}},
$$
and the right-hand side here belongs to $L^1(dA)$. Thus, the left-hand side of (3) is continuous, as well.

Then since $\|\cdot\|_{p,\alpha}$ is a norm, we have in view of (3) and (2) that
$$
\|\mathcal{H}_{K,A_\alpha}f\|_{p,\alpha}\le \int_{\mathbb{D}}|K(w)|\|f\circ \varphi_w\|_{p,\alpha}dA_\alpha(w)\le
$$
$$
\int_{\mathbb{D}}|K(w)|\left(\frac{1+|w|}{1-|w|}\right)^{\frac{2+\alpha}{p}}dA(w)\|f\|_{p,\alpha}.
$$
This proves the theorem.

\subsection{Hardy spaces}

As is well known, for $0 < p <\infty$ the \textit{Hardy space} $H^p=H^p(\mathbb{D})$ consists of analytic functions $f$ in the unit disc $\mathbb{D}$ such that
$$
\|f\|^{p}_{H^p}:= \sup_{0<r<1}\int_{0}^{2\pi} |f(re^{i\theta})|^pd\theta<\infty.
$$
Then $\|\cdot\|_{H^p}$ is a norm for $p\in [1,\infty)$.  It is easy to verify  also that this is a $p$-norm for $p\in (0,1)$, in particular for $\|\cdot\|_{H^p}^p$ the triangle inequality holds. Every  function $f(z)\in H^p(\mathbb{D})$ has  boundary values $f(e^{i\theta})\in L^p(\partial \mathbb{D})$ and
the map $f(z)\mapsto f(e^{i\theta})$ is an isometrical isomorphism of $H^p(\mathbb{D})$ onto some closed subspace of $L^p(\partial \mathbb{D},\frac{d\theta}{2\pi})$ (see, e.~g., \cite{Dur}, \cite{Zhu}).

Denote by $\widetilde{\Lambda}_p$
the class of all holomorphic functions $g$ on $\mathbb{D}$ which are continuous on the closure of $\mathbb{D}$ and such that if $1/p\notin\{2,3,\dots\}$ and $n=\lfloor 1/p\rfloor$ the   derivative $g^{(n-1)}$ of the $2\pi$ periodic counterpart of the boundary
function $g(e^{i\theta})$ belongs to the Lipschitz class $\Lambda_{\{1/p\}}$ on $\mathbb{R}$ (here $\{1/p\}$ stands for the fractional part of $1/p$) and  if $1/p=n+1$, then $g^{(n-1)}\in \Lambda_{\ast}$ where $\Lambda_{\ast}$ denotes the class  of all continuous functions $\varphi$ on $\mathbb{R}$ such that there is a constant
$A$ with the property
$$
|\varphi(x + h) - 2\varphi(x) +\varphi(x- h)| < Ah
$$
for all $x$ and for all $h > 0$.

\textbf{Theorem 3.}
 \textit{Let  the function  $K(w)/(1-|w|)^{1/p}$ belongs to $L^1(\mu)$. }

1) \textit{If $1\le p<\infty$, then the operator $\mathcal{H}_{K,\mu}$ is bounded on $H^p(\mathbb{D})$ and }
$$
\|\mathcal{H}_{K,\mu}\|_{H^p\to H^p}\le \int_{\mathbb{D}} |K(w)|\left(\frac{1+|w|}{1-|w|}\right)^{1/p}d\mu(w).
$$

2) \textit{Let $0<p<1$. The operator $\mathcal{H}_{K,\mu}$ is bounded on $H^p(\mathbb{D})$
if and only if for every $g\in \widetilde{\Lambda}_p$}
$$
c_n=\int_{\mathbb{D}}K(w)\left(\frac{g(0)}{\overline{w}^n}+
\underset{z=\overline{w}}{\mathrm{res}}\left(\left(\frac{wz-1}{z-\overline{w}}\right)^n\frac{g(z)}{z}\right)\right)d\mu(w),\quad n\in \mathbb{Z}_+\eqno(4)
$$
\textit{is  a sequens of Fourier coefficients of some function $h\in \widetilde{\Lambda}_p$.
}

Proof. 1) The inequality
$$
|f(z)|\le 2^{\frac{1}{p}}\|f\|_{H^p}(1-|z|)^{-\frac{1}{p}}\quad (f\in H^p(\mathbb{D}), z\in \mathbb{D})
$$
(see \cite[p. 36]{Dur}) shows that  the condition (a) of Lemma 1 is satisfied for $(X,m)=(\mathbb{D},dA)$ and $(\Omega,\mu)=(\mathbb{D},\mu)$. Conditions (b) and (c) are
the consequences of the inequality

$$
\|f\circ \varphi_w\|_{H^p}\le \left(\frac{1+|w|}{1-|w|}\right)^{1/p} \|f\|_{H^p}.
$$

 This  inequality follows from the Littlewood's subordination
theorem (see, e.g., \cite[Theorem 11.12]{Zhu}), but we shall give a simple direct proof.
If we as usual  identify the function $f(z)\in H^p$ with its boundary value $f(e^{i\theta})$, then
$$
\|f\circ \varphi_w\|_{H^p}^p=\frac{1}{2\pi}\int_{0}^{2\pi}|
f(\varphi_w(e^{i\theta}))|^pd\theta=\frac{1}{2\pi i}\int_{0}^{2\pi}|
f(\varphi_w(e^{i\theta}))|^p\frac{de^{i\theta}}{e^{i\theta}}.
$$
If we put in the last integral $e^{it}= \varphi_w(e^{i\theta})$, then $e^{i\theta}= \varphi_w(e^{it})$, and $de^{i\theta}= \frac{|w|^2-1}{(1-\overline{w}e^{it})^2}ie^{it}dt$.
Thus,
$$
\|f\circ \varphi_w\|_{H^p}^p=\frac{1-|w|^2}{2\pi}\int_{0}^{2\pi}\frac{|f(e^{it})|^pdt}{(\overline{w}-e^{-it})(w-e^{it})}=
$$
$$
(1-|w|^2)\frac{1}{2\pi}\int_{0}^{2\pi}\frac{|f(e^{it})|^pdt}{1+|w|^2-2\mathrm{Re}(we^{-it})}\le \frac{1+|w|}{1-|w|} \|f\|_{H^p}^p.
$$

 So, by Lemma 1 for $dA$ a.e. $z\in \mathbb{D}$
$$
\mathcal{H}_{K,\mu} f(z)=(B)\int_{\mathbb{D}}K(w)f\circ \varphi_wd\mu(w)(z),
$$
the Bochner integral for $H^p(\mathbb{D})$. As in the proof of Theorem 1 to show that this equality holds for all $z\in \mathbb{D}$ it suffices to prove the continuity of its left-hand side. To this end first note that by \cite[Theorem 9.1]{Zhu}
$$
|f(\varphi_w(z))|\le\frac{ \|f\|_{H^p}}{(1-|\varphi_w(z)|^2)^{1/p}}.
$$
On the other hand (see the proof of Theorem 2),
$$
1-|\varphi_w(z)|^2\ge \frac{1}{2}(1-|z|^2)(1-|w|).
$$
 Thus,
$$
|f(\varphi_w(z))|\le\frac{ 2^{1/p}\|f\|_{H^p}}{(1-|z|^2)^{1/p}(1-|w|)^{1/p}}.
$$
It follows, that for every compact neighborhood $U\subset \mathbb{D}$
$$
|K(w)||f(\varphi_w(z))|\le C_U\frac{|K(w)|}{(1-|w|)^{1/p}}
$$
is a $\mu$ integrable majorant for the left-hand side.
This proves the desired  continuity.

Since $\|\cdot\|_{H^p}$ is a norm, we have in view of the preceding  inequality that
$$
\|\mathcal{H}_{K,\mu} f\|_{H^p}\le \int_{\mathbb{D}}|K(w)|\|f\circ \varphi_w\|_{H^p}d\mu(w)\le
$$
$$
\int_{\mathbb{D}} |K(w)|\left(\frac{1+|w|}{1-|w|}\right)^{1/p}d\mu(w)\|f\|_{H^p}.
$$
This completes the proof of the first statement.

2) Let $0<p<1$. Note that $K\in L^1(\mu)$ Since $H^p(\mathbb{D})$ is an $F$-space and the dual $H^p(\mathbb{D})^*$ separates the points of $H^p(\mathbb{D})$ \cite[p. 118]{Dur}, Theorem II.2.7 in \cite{DS} shows that $\mathcal{H}_{K,\mu}$ is bounded on $H^p(\mathbb{D})$ if and only if a linear functional $l\circ \mathcal{H}_{K,\mu}$ belongs to $H^p(\mathbb{D})^*$ for each $l\in H^p(\mathbb{D})^*$.

We shall employ   the general form of a  linear
functional  on $H^p(\mathbb{D})$ \cite[Theorem 7.5]{Dur}.  By this theorem every bounded linear
functional  on $H^p(\mathbb{D})$ has a unique representation of the form
$$
l_g(f)=\lim_{r\to 1-0}\frac{1}{2\pi}\int_{0}^{2\pi} f(re^{i\theta})g(e^{-i\theta})d\theta,
$$
where $g\in \widetilde{\Lambda}_p$ (and vice versa).

Thus, $\mathcal{H}_{K,\mu}$ is bounded on $H^p(\mathbb{D})$ if and only if for every $g\in \widetilde{\Lambda}_p$ there is such $h\in \widetilde{\Lambda}_p$ that
$$
l_g\circ \mathcal{H}_{K,\mu}=l_h.\eqno(5)
$$
In other wards, for every $f\in H^p(\mathbb{D})$
$$
\lim_{r\to 1-0}\frac{1}{2\pi}\int_{0}^{2\pi}\left(\int_{\mathbb{D}}K(w)f(\varphi_w(re^{i\theta}))d\mu(w)\right)g(e^{-i\theta})d\theta
$$
$$ =\lim_{r\to 1-0}\frac{1}{2\pi}\int_{0}^{2\pi} f(re^{i\theta})h(e^{-i\theta})d\theta. \eqno(6)
$$

Let $\mathcal{H}_{K,\mu}$ is bounded on $H^p(\mathbb{D})$. Putting $f(z)=z^n$, $n\in \mathbb{Z}_+$ in  (6) we get in view  of the Fubini theorem that
$$
\lim_{r\to 1-0}\int_{\mathbb{D}}K(w)\left(\frac{1}{2\pi}\int_{0}^{2\pi}(\varphi_w(re^{i\theta}))^ng(e^{-i\theta})d\theta\right)d\mu(w)
$$
$$ =\lim_{r\to 1-0}\frac{1}{2\pi}\int_{0}^{2\pi} (re^{i\theta})^nh(e^{-i\theta})d\theta, \eqno(7)
$$
or by the Lebesgue theorem ($\varphi_w$ and $g$ are bounded)
$$
\int_{\mathbb{D}}K(w)\frac{1}{2\pi}\int_{0}^{2\pi}((\varphi_w(e^{i\theta}))^ng(e^{-i\theta})d\theta d\mu(w)
$$
$$=\frac{1}{2\pi}\int_{0}^{2\pi} e^{in\theta}h(e^{-i\theta})d\theta, \quad n\in \mathbb{Z}_+.\eqno(8)
$$
Note that by the Cauchy theorem
$$
\frac{1}{2\pi}\int_{0}^{2\pi}\left((\varphi_w(e^{i\theta})\right)^ng(e^{-i\theta})d\theta
=\frac{1}{2\pi}\int_{0}^{2\pi}\left(\frac{we^{it}-1}{e^{it}-\overline{w}}\right)^ng(e^{it})dt
$$
$$
 =\frac{1}{2\pi i}\int_{\{|z|=1\}}\left(\frac{wz-1}{z-\overline{w}}\right)^n\frac{g(z)}{z}dz
 =\frac{g(0)}{\overline{w}^n}+
 \underset{z=\overline{w}}{\mathrm{res}}\left(\left(\frac{wz-1}{z-\overline{w}}\right)^n\frac{g(z)}{z}\right).\eqno(9)
$$
Therefore (8) is equivalent to
$$
\int_{\mathbb{D}}K(w)\left(\frac{g(0)}{\overline{w}^n}+
\underset{z=\overline{w}}{\mathrm{res}}\left(\left(\frac{wz-1}{z-\overline{w}}\right)^n\frac{g(z)}{z}\right)\right)d\mu(w)
$$
$$
=
\frac{1}{2\pi}\int_{0}^{2\pi} e^{-in\theta}h(e^{i\theta})d\theta, \quad n\in \mathbb{Z}_+ \eqno(10)
$$
and the necessity follows.

To prove the sufficiency, for every $g\in \widetilde{\Lambda}_p$ let $h\in \widetilde{\Lambda}_p$ be such that for  its Fourier coefficients $c_n$ (4) holds.
Then formula (10) is valid. In view of (9) this implies (8). In turn, (8) implies (7) and then for every
algebraic polynomial $q_n$ we get
$$
\int_{\mathbb{D}}K(w)\left(\lim_{r\to 1-0}\frac{1}{2\pi}\int_{0}^{2\pi}q_n(\varphi_w(re^{i\theta}))g(e^{-i\theta})d\theta\right)d\mu(w)
$$
$$
 =\lim_{r\to 1-0}\frac{1}{2\pi}\int_{0}^{2\pi} q_n(re^{i\theta})h(e^{-i\theta})d\theta.
$$
 In other words,
$$
\int_{\mathbb{D}}K(w)l_g(q_n\circ \varphi_w)d\mu(w)=l_h(q_n).\eqno(11)
$$
 Taking into account that polynomials are dense in $H^p$ (see, e.~g., \cite[Corollary 9.5]{Zhu}), for every $f\in H^p$ one  can choose   a sequence $q_n$ of polynomials
that converges to $f$ in $H^p$. Since
$$
l_h(q_n)\to l_h(f), \ l_g(q_n\circ \varphi_w)\to l_g(f\circ \varphi_w) \mbox{ as } n\to\infty,
$$
 and by Lemma 1
$$
|l_g(q_n\circ \varphi_w)|\le \|l_g\|\|q_n\|_p\left(\frac{1+|w|}{1-|w|}\right)^{1/p}\le \mathrm{const}\left(\frac{1+|w|}{1-|w|}\right)^{1/p},
$$
formula (11) implies in view of the Lebesgue theorem that the property (5) is valid. This proves the sufficiency.$\Box$

Putting $g(z)=z$ in the previous theorem we get the following

\textbf{Corollary 1.} \textit{Let $0<p<1$ and let the function  $K(w)/(1-|w|)^{1/p}$ belongs to $L^1(\mu)$. If the operator $\mathcal{H}_{K,\mu}$ is bounded on $H^p(\mathbb{D})$, then }
$$
c_n=n\int_{\mathbb{D}}K(w)w^{n-1}(|w|^2-1)d\mu(w),\quad n\in \mathbb{Z}_+
$$
\textit{is  a sequens of Fourier coefficients of some function $h\in \widetilde{\Lambda}_p$.
}

\section{Identity Approximation by Hausdorff Operators }

Although the class of operators we are considering differs  from the class considered in \cite{KSZhu}, there are analogs of Theorem 14 from \cite{KSZhu} that was devoted to identity approximation by Hausdorff–Berezin operators.

In the following for the function $f$ on $\mathbb{D}$ we put $f^\vartriangle(z)=f(-z)$ and
for $0<\varepsilon <1$  consider  operators  of the form
$$
\mathcal{H}_\varepsilon f(z):=\int_{\mathbb{D}}K(w)(f\circ \varphi_{\varepsilon w})^\vartriangle(z) d\mu(w). \eqno(12)
$$

As in \cite[Section 5]{KSZhu} for the goals of approximation it is naturally to assume that
$$
K\in L^1(d\mu),\ \mbox{ and } \ \int_{\mathbb{D}}Kd\mu=1. \eqno (13)
$$

\textbf{Theorem 4}. \textit{Let $1\le p < \infty$. Under the assumption (13) with $dA_\alpha$ instead  of $d\mu$ the operators (12) with $dA_\alpha$ instead  of $d\mu$ are
identity approximations in $L^p(dA_\alpha)$, namely,}
$$
\lim\limits_{\varepsilon \to 0}\|\mathcal{H}_\varepsilon f-f\|_{p,\alpha}=0. \eqno(14)
$$

Proof. As in the proof of Theorem 2 we have
$$
\mathcal{H}_\varepsilon f=(B)\int_{\mathbb{D}}K(w)(f\circ \varphi_{\varepsilon w})^\vartriangle dA_\alpha(w), \eqno(15)
$$
the Bochner integral in  $L^p(dA_\alpha)$. It follows in view of (13) that
$$
\|\mathcal{H}_\varepsilon f-f\|_{p,\alpha}\le \int_{\mathbb{D}}|K(w)|\|(f\circ \varphi_{\varepsilon w})^\vartriangle -f \|_{p,\alpha}dA_\alpha(w). \eqno(16)
$$
 We shall show that the operators $\mathcal{H}_\varepsilon$ are uniformly bounded when $\varepsilon\in(0,1/2)$, specifically,
$$
\|\mathcal{H}_\varepsilon\|\le 2^{\frac{2(2+\alpha)}{p}}\int_{\mathbb{D}}|K(w)|dA_\alpha(w). \eqno(17)
$$
Indeed, using
  \cite[Proposition 4.3]{Zhu} with $a=\varepsilon w$  we have
$$
\|(f\circ \varphi_{\varepsilon w})^\vartriangle)\|_{p,\alpha}^p=\|f\circ \varphi_{\varepsilon w}\|_{p,\alpha}^p=
$$
$$
\int_{\mathbb{D}}|f(z)|^p\frac{(1-|\varepsilon w|^2)^{(2+\alpha)}}{|1-\varepsilon \overline{w} z|^{2(2+\alpha)}}dA_\alpha(z)\le
$$
$$
\int_{\mathbb{D}}|f(z)|^p\frac{1}{(1-\varepsilon)^{2(2+\alpha)}}dA_\alpha(z)\le 2^{2(2+\alpha)}\|f\|_{p,\alpha}^p.\eqno(18)
$$
Thus, (17) follows from (15) and (18).

In view of (17) to prove (14)  it remains to check this equality on a dense subset
of $L^p(dA_\alpha)$. We shall use  (16) to check (14) on the subset of bounded functions.
Let $|f(z)|\le C$. Then
$$
\|(f\circ \varphi_{\varepsilon w})^\vartriangle -f \|_{p,\alpha}^p=\int_{\mathbb{D}}|(f\circ\varphi_{\varepsilon w})(-z)-f(z)|^pdA_\alpha(z)\to 0
$$
as $\varepsilon\to 0$ by the Lebesgue's dominated convergence theorem ($|(f\circ\varphi_{\varepsilon w})(-z)-f(z)|\le 2C$). But formula (18) implies that
$$
\|(f\circ \varphi_{\varepsilon w})^\vartriangle -f \|_{p,\alpha}\le (2^{\frac{2(2+\alpha)}{p}}+1) \|f\|_{p,\alpha}.
$$
It follows  (again by the Lebesgue's dominated convergence theorem) that the right-hand side of (16) vanishes as $\varepsilon\to 0$. This completes the proof.

\textbf{Theorem 5}. \textit{Let $1\le p < \infty$. Under the assumption (13)  the operators (12) are
identity approximations in $H^p(\mathbb{D})$, namely,}
$$
\lim\limits_{\varepsilon \to 0}\|\mathcal{H}_\varepsilon f-f\|_{H^p}=0. \eqno(19)
$$

The proof is similar with the proof of Theorem 4.  Indeed, as in the proof of Theorem 3 we have
$$
\mathcal{H}_\varepsilon f=(B)\int_{\mathbb{D}}K(w)(f\circ \varphi_{\varepsilon w})^\vartriangle d\mu(w), \eqno(20)
$$
the Bochner integral in  $H^p(\mathbb{D})$. It follows in view of (13) that
$$
\|\mathcal{H}_\varepsilon f-f\|_{H^p}\le \int_{\mathbb{D}}|K(w)|\|(f\circ \varphi_{\varepsilon w})^\vartriangle -f \|_{H^p}d\mu(w). \eqno(21)
$$
 We shall show that the operators $\mathcal{H}_\varepsilon$ are uniformly bounded when $\varepsilon\in(0,1/3)$, specifically,
$$
\|\mathcal{H}_\varepsilon\|\le 2^{\frac{1}{p}}\int_{\mathbb{D}}|K(w)|d\mu(w). \eqno(22)
$$
Indeed, in the proof of Theorem 3 it was shown that
$$
\|(f\circ \varphi_{\varepsilon w})^\vartriangle)\|_{H^p}=\left(\frac{1+\varepsilon |w|}{1-\varepsilon |w|}\right)^{\frac{1}{p}}\|f\|_{H^p}.
$$
Thus,
$$
\|(f\circ \varphi_{\varepsilon w})^\vartriangle)\|_{H^p}\le 2^{\frac{1}{p}}\|f\|_{H^p} \eqno(23)
$$
 as $\varepsilon\in(0,1/3)$,  and  (22) follows from (20) and (23).

In view of (22) to prove (19)  it remains to check this equality on a dense subset
of $H^p(\mathbb{D})$. We shall use  (21) to check (19) on the subset of polynomials.
First note that since $|z_1^k-z_2^k|\le (k+1)|z_1-z_2|$ for $z_1, z_2\in \mathbb{D}$, $k\in \mathbb{N}$, we have
$$
|((\varphi_{\varepsilon w}(z))^k)^\vartriangle)-z^k|=\left|\left(\frac{\varepsilon w+z}{1+\varepsilon \overline{w}z}\right)^k-z^k\right|\le (k+1)\left|\frac{\varepsilon w-\varepsilon \overline{w}z^2}{1+\varepsilon \overline{w}z}\right|\le \frac{2(k+1)\varepsilon}{1-\varepsilon}.
$$
and therefore
$$
\|((\varphi_{\varepsilon w}(z))^k)^\vartriangle)-z^k\|_{H^p}\le \frac{2(k+1)\varepsilon}{1-\varepsilon}.
$$
It follows that
$$
\|(f\circ \varphi_{\varepsilon w})^\vartriangle -f \|_{H^p}\to 0
$$
as $\varepsilon\to 0$ for every polynomial $f$. But formula (23) implies that
$$
\|(f\circ \varphi_{\varepsilon w})^\vartriangle -f \|_{H^p}\le (2^{\frac{1}{p}}+1) \|f\|_{H^p}.
$$
 Then the Lebesgue's dominated convergence theorem shows that the right-hand side of (21) vanishes as $\varepsilon\to 0$. This completes the proof.



\end{document}